\documentclass[a4paper,12pt]{article}%,openany
\usepackage{amsfonts,amssymb,amsmath,amsthm}%,naglowki}
\usepackage[dvips]{graphics}
\usepackage[dvips]{graphicx}
\usepackage[dvips]{color}

\def\num{\hspace{-2mm}{\bf }\hspace{2mm}}

\newtheorem{st}{Statement}[section]

\newtheorem{prob}{Problem}
\newtheorem{propo}[st]{Proposition}

\newtheorem{cor}[st]{Corollary}
\newtheorem{rem}[st]{Remark}
\newtheorem{thm}[st]{Theorem}
\newtheorem{lemm}[st]{Lemma}

\def\max{{\rm max\,}}
\def\min{{\rm min\,}}

\def\conf{{\looparrowleft}}
\def\cP{{\cal P}}

\def\card{{\rm card\,}}

\def\Proof:{ \vspace{-1.5mm} {\noindent\it Proof.}}
\def\gm{\vspace{0mm}}%-2mm
\def\gd{\vspace{0mm}}%-5mm

\def\Box{\rule{1.5mm}{1.5mm}}

\begin{document}

\renewcommand{\thefootnote}{\fnsymbol{footnote}}

 \title{\gd \gd Ordinal ultrafilters versus P-hierarchy }
 \author{ Andrzej Starosolski}
 \date{\today}
\maketitle

\gd
\begin{abstract}

An earlier paper, entitled ``P-hierarchy on $\beta\omega$'',
investigated the relations between ordinal ultrafilters and the
so-called P-hierarchy. 
The present
paper
focuses on the aspects of characterization of classes of ultrafilters of finite
index, existence, generic existence and the Rudin-Keisler-order.

\end{abstract}


\gd
\section{Introduction}

Ultrafilters on $\omega$ may be classified with respect to
sequential contours of different ranks, that is, iterations of the
Fr\'{e}chet filter by contour operations. This way an
$\omega_1$-sequence $\{\cP_\alpha \}_{1\leq\alpha\leq\omega_1}$ of
pairwise disjoint classes of ultrafilters - the P-hierarchy - is
obtained, where P-points correspond to the class $\cP_2$, allowing
us to look at the P-hierarchy as on extension of P-points. Section
2 recalls all necessary definitions and properties of the
P-hierarchy. Section 3 shows some equivalent conditions for an
ultrafilter to belong to a class of (fixed) finite index of the
P-hierarchy; those conditions appear to be very similar to the
behavior of classical P-points. We also obtain another condition for
belonging to a class of (fixed) finite index of the P-hierarchy
which is literally a part of conditions for being an element of a
class of (fixed) finite index of ordinal ultrafilters. Section 4
focuses on the Rudin-Keisler order on P-hierarchy classes. It is
shown that RK minimal elements of classes of finite index can exist.
Similar results are achieved for ordinal
ultrafilters. In section 5 we show evidence for the generic existence
of the P-hierarchy being equivalent to $\mathfrak{d}=\mathfrak{c}$,
in consequence, being equivalent to the generic existence of ordinal
ultrafilters. In section 6 we prove that CH implies that each class
of the P-hierarchy is not empty, we also presented known results concerning
existence of both types of ultrafilters.

We generally use standard terminology, however less popular terms are taken from \cite{DolMyn},
where key-term "monotone sequential cascades" has been introduced. All necessary information
may also be found in \cite{Star-P-hier}.
For additional information regarding sequential cascades and
contours a look at \cite{DolMyn}, \cite{DolStaWat},
\cite{Dol-multi}, \cite{Star-ff} is recommended. Below, only the
most important definitions and conventions are repeated.

If $u$ is a filter on $A \subset B$, then we identify $u$ with the
filter on $B$ for which $u$ is a filter-base.

Let $p$ be a filter on $X$, and let $q$ be a filter on
$Y$; we say that $p$ is {\it Rudin-Keisler} greater than
$q$ (we write $p \geq_{RK} q$) if there is such a
map $f: X \rightarrow Y$ that $f(p) \supset q$. We say
that $p$ is {\it infinite Rudin-Keisler} greater than $q$ (we write $p >_{\infty} q$) if there is a map $f:X
\rightarrow Y$ with $f(p) =  q$, but
%and
there is no
$P \in p$ such that $f\mid_P$ is finite-to-one. We say that
$p$ is {\it greater} than $q$ if
$q \subset p$.

Recall also that if $p$, $q$ are ultrafilters, $f(p) = q$ and if $p \approx q$ (i.e. $p
\geq_{RK} q$ and $q \geq_{RK} p$), then there
exists $P \in p$ such that $f\mid_p$ is one-to-one (see
\cite[Theorem 9.2]{ComNeg}).

The {\it cascade} is a well founded tree i.e. a tree $V$ without infinite branches and with a
least element $\emptyset _V$. A cascade is $\it sequential$ if for
each non-maximal element of $V$ ($v \in V \setminus \max V$) the set
$v^{+V}$ of immediate successors of $v$ (in $V$) is countably
infinite. For $v \in V$ we write $v^{-V}$ to denote such an element of $V$ that $v\in (v^{-V})^{+V}$.
For $A\subset V$ we use $A^{+V} = \bigcup_{v\in A}v^{+V}$, $A^{-V} = \bigcup_{v\in A}v^{-V}$.
In symbols $v^{+V}$, $v^{-V}$, $A^{+V}$, $A^{-V}$ we omit the name of cascade (obtaining $v^{+}$, $v^{-}$, $A^{+}$, $A^{-}$)if it is clear from the context which cascade we have on mind. 
If $v \in V \setminus
\max V$, then the set $v^+$ (if infinite) may be endowed with an
order of the type $\omega$, and then by $(v_n)_{n \in \omega}$ we
denote the sequence of elements of $v^+$, and by $v_{nV}$ - the
$n$-th element of $v^{+V}$.

The {\it rank} of $v \in V$ ($r_V(v)$ or $r(v)$) is defined
inductively as follows: $r(v)=0$ if $v \in \max V$, and otherwise
$r(v)$ is the least ordinal greater than the ranks of all immediate
successors of $v$. The rank $r(V)$ of the cascade $V$ is, by
definition, the rank of $\emptyset_V$. If it is possible to order
all sets $v^+$ (for $v \in V \setminus \max V$)  so that for each $v
\in V \setminus \max V$ the sequence $(r(v_n)_{n<\omega})$ is
non-decreasing, then the cascade $V$ is {\it monotone}, and we fix
such an order on $V$ without indication.

Let $W$ be a cascade, and let $(V_w)_{w \in \max W}$ be a pairwise
disjoint sequence of cascades such that $V_w \cap W = \emptyset$ for
all $w \in \max W$. Then, the {\it confluence} of cascades $V_w$
with respect to the cascade $W$ (we write $W \looparrowleft V_w$) is
defined as a cascade constructed by the identification $w \in \max
W$ with $\emptyset_{V_w}$ and according to the following rules:
$\emptyset_W = \emptyset_{W \looparrowleft V_w}$; if $w\in W
\setminus \max W$, then $w^{+ W \looparrowleft V_w} = w^{+W}$; if $w
\in V_{w_0}$ (for a certain $w_0 \in \max W$), then $w^{+ W
\looparrowleft V_w} = w^{+V_{w_0}}$; in each case we also assume
that the order on the set of successors remains unchanged. By $(n)
\looparrowleft V_n$ we denote $W \looparrowleft V_w$ if $W$ is a
sequential cascade of rank 1.

If $\mathbb{P}= \{p_s: s \in S \}$ is a family of filters on
$X$ and if $q$ is a filter on $S$, then the {\it contour of $\{
p_s \}$ along $q$} is a filter on $X$ defined by
$$\int_q \mathbb{P} = \int_{q}p_s =
\bigcup_{Q \in q} \bigcap_{s \in Q} p_s.$$

Such a construction has been used by many authors (\cite{Fro},
\cite{Gri1}, \cite{Gri2}) and is also known as a sum (or as a limit)
of filters. On the sequential cascade, we consider the finest
topology such that for all but the maximal elements $v$ of $V$, the
co-finite filter on the set $v^{+V}$ converges to $v$. For the
sequential cascade $V$ we define the {\it contour} of $V$ (we write
$\int V$) as the trace on $\max V$ of the neighborhood filter of
$\emptyset _V$ (the {\it trace} of a filter $u$ on a set $A$ is the
family of intersections of elements of $u$ with $A$). Equivalently we may sey
that $\int V$ is a Fr\'echet filter on $\max V$ if $r(V)=1$, and $\int V = \int_{Fr}\int V_n$ if 
$V= (n) \looparrowleft V_n$ where $Fr$ denotes the Fr\'echet filter.  Similar
filters were considered in \cite{Kat1}, \cite{Kat2}, \cite{Dagu}.
Let $V$ be a monotone sequential cascade and let $u=\int V$.
Then
the {\it rank $r(u)$} of $u$
%($r(u)$)
is, by definition, the rank of $V$.
%, it
It
was shown in \cite{DolStaWat} that if $\int V= \int W$, then $r(V) =
r(W)$.
%).

Let $S$ be a countable set. A family $\{u_s\}_{s\in S}$ of filters
is referred to as $\it discrete$ if there exists a pairwise disjoint
family $\{U_s\}_{s \in S}$ of sets such that $U_s \in u_s$ for each
$s \in S$. For $v \in V$ we denote by $v^\uparrow$ a subcascade of
$V$ built by $v$ and all successors of $v$. If $U \subset \max V$
and $U \in \int V$, then by $U^{\downarrow V}$ we denote the biggest
(in the set-theoretical order)$^{1)}$ {\footnotetext{ i.e. the cascade $V$
is greater then the cascade $W$ if $V \subset W$, $\emptyset_V=\emptyset_W$,
$w^{+W}\subset w^{+V}$ for all $w \in W$, $\max W \subset \max V$.}}
monotone sequential subcascade of
the cascade $V$ built of some $v \in V$ such that $U \cap \max v
^\uparrow \not= \emptyset$.  We write $v^{\uparrow }$ and
$U^\downarrow$ instead of $v^{\uparrow V}$ and $U^{\downarrow V}$ if
we know in which cascade the subcascade is considered. By $V_n$ we usually denote
$(\emptyset_V)_n^\uparrow$, by $V_{n,m}$ we understand  $((\emptyset_V)_n^+)_m^\uparrow$.

\section{Ordinal ultrafilters and classes $\cP_\alpha$}

In the remainder of this paper each filter is considered to be on
$\omega$, unless otherwise indicated. Let us define ${\cP}_\alpha$
for $1 \leq {\alpha <\omega_1}$ on $\beta\omega$ (see
\cite{Star-P-hier}) as follows: $u \in {\cP}_\alpha$ if there is no
monotone sequential contour ${\cal V}_{\alpha}$ of rank $\alpha$ such that
${\cal V}_{\alpha}  \subset u$, and for each $\beta$ in the range
$1\leq\beta < \alpha$ there exists a monotone sequential contour
${\cal V}_\beta$ of rank $\beta$ such that ${\cal V}_\beta \subset u$. Moreover,
if for each $\alpha < \omega_1$ there exists a monotone sequential
contour ${\cal V}_\alpha$ of rank $\alpha$ such that ${\cal V}_\alpha \subset u$,
then we write $u \in {\cP}_{\omega_1}$.

Let us recall three equivalent definitions of P-points:  a point $u
\in \beta\omega \setminus \omega$ is a {\it P-point} if \newline A)
the intersection of countably many neighborhoods of $u$ is a (not
necessarily open) neighborhood of $u$; \newline B) for each
countable set $\{U_n\}_{n<\omega}$  of elements of the ultrafilter
$u$ there exists a set $U \in u$ such that $\card(U \setminus U_n)$
is finite for each $n<\omega$;
\newline C) for each function $f:\omega \rightarrow \omega$ there
exists a set $U \in u$ such that either $f\mid_U$ is constant or
$f\mid_U$ is finite-to-one.

\gm\begin{rem}\num If $u$ is an ultrafilter on $\omega$ then:

1) $u \in {\cP}_1$ if and only if $u$ is a principal ultrafilter;

2) if $u$ is RK-minimal then $u \in {\cP}_2$. $\Box$
\end{rem}

%For $\alpha < \omega_1$, let us define an {\it $\alpha$-partition}
%as follows: let
Let $M$ be a countably infinite set, and let $V$ be a monotone
cascade of rank $\alpha<\omega_1$ such that $\max V = M$.
%The
Then the
set $D = \{ D_v =  \max v^\uparrow : v \in V, r(v) \geq 1 \}$ is
called an {\it $\alpha$-partition} (of $M$).

Thus, the classic ``partitions of $\omega$ into infinitely many
infinite sets'' belong to ``2-partitions'' in our language. Since a
cascade uniquely defines a partition, it is usually identified with
its cascade. For an $\alpha$-partition we define by transfinite
induction {\it residual sets} as follows: a set $A$ is residual for
a 1-partition $V$ if $A \cap \max V$ is finite; if residual sets are
defined for all $\beta$-partitions for $\beta < \alpha$, then a set
$A$ is residual for the $\alpha$-partition $V=(n) \looparrowleft V_n$ 
 if there exists a finite set $N \subset \omega$ such that for all
$n \not\in N$ the set $A$ is residual for the
partitions $V_n$. For a partition defined by a monotone sequential
cascade $V$, equivalently we can say that $U$ is residual if and
only if $\omega\setminus U \in \int V$.

Certain properties of the P-hierarchy from \cite{Star-P-hier} are
listed below, namely Proposition 2.1 and Theorems 2.3, 2.5, 2.9,
2.8.

\begin{propo}\num An ultrafilter $u$ (on $\omega$) is a P-point if and
only if $u \in {\cP}_2$.
\end{propo}

\gm\begin{thm}\num Let $u \in {\cP}_\alpha$ and let $f:\omega
\rightarrow  \omega$. Then $f(u) \in {\cP}_\beta$ for a certain
$\beta \leq \alpha$.
\end{thm}

Let $\alpha$ be an ordinal, by $-1+\alpha$ we understand $\alpha -
1$ if $\alpha$ is finite, and $\alpha$ if $\alpha$ is infinite.

\begin{thm}\num Let $(\alpha_n)_{n<\omega}$ be a non-decreasing
sequence of ordinals less than $\omega_1$, let $\alpha =
\lim\limits_{n < \omega} (\alpha_n )$, let $1 < \beta <\omega_1$ and
let $(X_n)$ be a partition of $\omega$. If $(p_n)$ is a sequence of
ultrafilters such that $X_n \in p_n \in \cP_{\alpha_n}$ and $p \in
\cP_{\beta}$, then $\int_p p_n \in \cP_{\alpha+(-1+\beta)}$.
\end{thm}

\begin{thm}\num
Let $\alpha$, $\beta$ be countable ordinals. If $u \in
\cP_{\alpha+\beta+1}$ then there exists a function $f:\omega
\rightarrow \omega$ such that $f(u) \in \cP_{\beta+1}$.
\end{thm}

\begin{thm}\num The following statements are equivalent:

1) P-points exist,

2) the $\cP_\alpha$ classes are non-empty for each countable successor
$\alpha$,

3) There exists a countable  successor $\alpha>1$ such that the
class $\cP_\alpha$ is non-empty.
\end{thm}

In \cite{Baum} Baumgartner provides the following definition.
%:
Let
$I$ be a family of subsets of a set $A$ such that $I$ contains all
singletons and is closed under subsets. Given an ultrafilter $u$ on
$\omega$, we say that $u$ is an $I$-ultrafilter if for any $f:
\omega \rightarrow A$ there is $U \in u$ such that $f(U) \in I$. For
$\alpha < \omega_1$, let $I_\alpha= \{ B \subset \omega_1: B$ has an
order type $\leq \alpha \}$, $J_\alpha = \{ B \subset \omega_1: B$
has order type $< \alpha \}$. %For convenience, we refer to
A {\it proper} $I_\alpha$-ultrafilter is one which is not an
$I_\beta$-ultrafilter for any $\beta < \alpha$. Denote also proper
$J_\alpha$-ultrafilters as $J^*_\alpha$-ultrafilters those being the
$J_\alpha$-ultrafilters which are not $J_\beta$-ultrafilters for any
$\beta < \alpha$. We can also find in $\cite{Baum}$ the following
statement:
%that:
If $u$ is
a proper $I_\alpha$-ultrafilter, then $\alpha$ must be
indecomposable. Recall that

\begin{propo}\cite[a corollary of Proposition 3.3]{Star-P-hier} If $u \in
J^*_{\omega^\alpha}$, then $u \in \cP_\beta$ for a certain $\beta
\leq \alpha$.
\end{propo}

\section{$\cP_\alpha$ classes for finite $\alpha$ and $<_{\infty}$ sequences}

\begin{thm}\num
If $u \in \beta\omega$, then the following statements are
equivalent:

1) There is no monotone sequential contour $C$ of rank $n$ such that
$C \subset u$. (i.e., for each n-partition there exists a set $U \in
u$ residual for this partition)

2) $u \in \bigcup\limits_{i=1}^n \cP_i$.

3) For each family of functions $\{f_1,\ldots , f_{n-1}\}$,
$f_i:\omega\rightarrow\omega$ there exists a set $U \in u$ such that
\newline
a) $f_1 \circ \ldots \circ f_{n-1} \mid_U$ is constant
    or \newline
b) there exists $i \in \{1, \ldots, n-1\}$ such that
    $f_i\mid_{f_{i+1}\circ \ldots \circ f_{n-1}(U)}$ is finite-to-one.

4) For each function $f:\omega\rightarrow\omega$ there exists a set
$U \in u$ such that
\newline
a) $f^{n-1}\mid_U$ is constant
    or \newline
b) there exists $i \in \{1, \ldots, n-1\}$ such that
    $f\mid_{f^{i-1}(U)}$ is finite-to-one.
\end{thm}

\Proof: $1\Leftrightarrow  2$ is trivial.

$2\Rightarrow3$: Let $u \in \cP_i$ for some $i\leq n$ and let us take any functions
$f_1, \ldots, f_{n-1}: \omega \rightarrow \omega$.

Let $A_k^{\infty}= \{m < \omega: \card(f^{-1}_k(m))= \omega \}$
and
$A_k^{fin}= \{m < \omega: \card(f^{-1}_k(m))< \omega \}$ for $k \in
\{1, \ldots, n-1\}$. Since $u$ is an ultrafilter, and $A_k^{\infty}
\cup A_k^{fin} = \omega$, for each $k$ one of those sets belongs to  $f_{k+1} \circ \ldots \circ f_{n-1}(u)$.
If for some $k$ it is $A_k^{fin}$,
then case 3b) holds, so we can assume that for each $k$, each function
$f_k$ is infinite-to-one on elements of $f_{k+1}\circ \ldots \circ
f_{n-1} (u)$. Since our research
%are
is
restricted to elements of
images of $u$, without loss of generality we may assume that
$\card(f^{-1}_k(m))=\omega$ for each $k\in \{1, \ldots, n-1\}$ and
for each $m \in \omega$.

Note the following obvious claim: Let $u$ be an ultrafilter and let $f$ be a function such that $f^{-1}(n)$ is infinite for all $n<\omega$. Then for each monotone sequential cascade $V$ of rank $\alpha$ such that $\int V \subset f(u)$, there is $\int f^{-1}(V) \subset u$, and $r(f^{-1}(v))=1+\alpha$, where $f^{-1}(V)= V \conf f^{-1}(v)$. 

If $f_{1}\circ \ldots \circ f_{n-1}(u)$ is not a principal ultrafilter, then $f_{k}\circ \ldots \circ f_{n-1}(u)$ contain a contour of rank $k$, and thus $u$ contain a contour of rank $n$ - a contradiction.

$3\Rightarrow4$ is trivial.

$4\Rightarrow1$: Let us assume that there exists a monotone
sequential contour $C_{n}$ of rank $n$ such that $C_{n} \subset u$.
There exists a monotone sequential cascade $V$
% exists,
such that $\int V =
C_{n}$. Naturally, $r(V)=n$. Without loss of generality we may assume
that $\max V = \omega$ and the cascade $V$ is complete, i.e. each
branch has the same length $n$. We identify elements of $\max V$
with $n$-sequences of natural numbers which label that elements i.e. $\emptyset_V=\emptyset$, $({i_1,\ldots, i_k})^+_{i_{k+1}} = (i_1,\ldots, i_{k+1})$. We define the function $f:\max V\rightarrow \max V$
as follows:

$f((i,1,1,...,1,1))= (1,1,...,1,1)$ for each $i<\omega$;

if $v=(k_1,\ldots, k_{n})$ and if there exists $l\in \{2,\ldots,n\}$
such that $k_l\not= 1$ then let $m(v)=\min\{t\in\{1,\ldots,n\}:
\forall l>t, l\leq n$: $ k_l=1 \}$ and let $f(v)=
(k_1,\ldots,k_{m(v)-2},k_{m(v)-1}+1,1,\ldots,1,1)$.

It may be noticed without difficulty that $f^i(\max V)  = \{v \in \max V: m(v)<n-i\}$ for $i\in\{0,\ldots,n-1\}$. Let $V(i)=
(f^i(\max V))^{\downarrow V}$.

Take any $G \subset f^i(\max V)$. If $f$ is finite-to-one on $G$, then
$G\cap \max v^{\uparrow V(i)}$ is finite for each $v\in V(i)$ such
that $v$ is a sequence of length $n-i$. Thus, $G$ is residual for $V(i)$ and so does not
belong to $f^i(u)$.

Thus, the ultrafilter $u$ does not have the property described in
point 4 of Theorem 3.1. $\Box$

\hspace{3.5mm}

It is worth comparing the definitions of P-points from page 3 with
the conditions of Theorem 3.1 in order to see that the behavior of
P-points is, in a very natural way, extended onto the behavior of
elements of $\bigcup\limits_{i=1}^n \cP_i$. Condition 1 of Theorem
3.1 is the extension of the equivalent definition of P-point from
Theorem 2.3, Condition 2 can be expressed as ``$u$ is no more than
${\cal P}_n$-point'', and Conditions 3 and 4 of Theorem 3.1 extend the
definition ``C'' of P-point from page 3.

\begin{propo}\num Let $u \in \cP _n$, $n \in \omega$, and
$f:\omega\rightarrow\omega$. If $f(u) \in \cP_n$ then there exists a
set $U \in u$ such that $f\mid_U$ is finite-to-one.
\end{propo}
\Proof: Proof basis on the same idea as a proof of $2\Rightarrow 3$ in the previous Theorem 3.1.
Suppose on the contrary, that there is no $U\in u$
that $f\mid_U$ is finite-to-one, thus $\{ i<\omega: \card(f^{-1}(i))=\omega\}\in f(u)$, and so without loss of generality we may assume that $\{ i<\omega: \card(f^{-1}(i))=\omega\}=\omega$. 
If $f(u)\in \cP_n$, then there exists a monotone sequential
contour ${\cal V} \subset f(u)$ of rank $n-1$. Consider a monotone
sequential cascade $V$ such that $\int V = {\cal V}$ and
$W=V\looparrowleft f^{-1}(v)$, where $v \in \max V$.
Since $W$ is a cascade of rank $n$ and $u\in \cP_n$, there exists a
set $U\in u$ residual for $W$. 
Consider sets $W^i$ for $i\in \{0,
\ldots, n-1\}$ defined by $W_0=U$ and $w \in
W^{i+1}\Leftrightarrow\card (w^{+W} \cap W^i) = \omega$ (sets $W^i$
are subsets of $W$ and  of $V$ as well, for $i > 0$). Split $U$
into $n$ pieces: $U_{n-1} = (W^{n-1})^{\uparrow W} \cap U$, $U^{i-1}=
((W^{i-1})^{\uparrow W} \cap U) \setminus \bigcup_{j=i}^n U_k$. Notice
that $U_i \not\in u$ for $i>0$. Indeed, $f(U_i) \subset
(W^i)^{\uparrow V} \cap \max V$ and $(W^i)^{\uparrow V} \cap \max V$ is
residual for $V$. Thus $U_0\in u$, clearly $f\mid_{U_0}$ is
finite-to-one. $\Box$

\hspace{3mm}

By Proposition 3.2 and Theorems 2.3 and 2.5 we obtain the following

\begin{cor}
If $u$ is an ultrafilter, then $u \in \cP_n$ if and only if there
exists an $n$-element $<_\infty$-decreasing sequence below it that contain "u" and a principal ultrafilter, and
there is
no
such chain of length $n+1$.
\end{cor}

\Proof: Non existence of such $n+1$ chain follows from Proposition 3.2. Existence of $n$ chain follows inductively from a following fact:
If $\int V \subset u$ for monotone sequential cascade $V$, then $\int f(V) \subset f(u)$ where $f[v^+]=v_1$ for all $v\in V: r(v)=1$, and $f$ is an identity on the rest of $V$ - such defined $f$ is not finite-to-one on each $U\in u$ (for details see proof of Proposition 3.5). Note that if $r(V)$ is finite than $r(f(V))=r(V)-1$. $\Box$

In \cite{Laf} Laflamme shows:

\begin{propo}\cite[Reformulation of Proposition 2.3]{Laf}

Let $k\in \omega$ and $u$ an ultrafilter such that

(*) $(\forall h \in \, {} ^\omega\omega_1)(\exists X \in u)$ the
order type of $h(X)$ is strictly less than $\omega^\omega$

Then $u$ is an $J^*_{\omega^k}$-ultrafilter precisely if it has a
$<_\infty$-chain of length $k$ below it that contain $u$ and a principal ultrafilter but
no such a chain of length $k+1$.
\end{propo}

Notice that Proposition 3.4 for ordinal ultrafilters is very similar
to Corollary 3.3, the only difference being the extra assumption (*).

As opposed to Proposition 3.2, for infinite $\alpha$'s we have the
following

\begin{propo}\num
If  $\alpha$ is a countably infinite successor ordinal, then for
each $u \in \cP_\alpha$ there exists a function $f:\omega\rightarrow
\omega$ such that: $f(u) \in \cP_\alpha$ and $f\mid_U$ is not
finite-to-one for each $U\in u$.
\end{propo}
\Proof: Let $u \in \cP_\alpha$, let $\alpha$ be as in the
assumptions. Let us take a monotone sequential contour
${\cal V}$ of rank $\alpha-1$ such that ${\cal V} \subset u$;
consider a monotone sequential cascade $V$ such that $\int V =
{\cal V}$; without loss of generality we may assume that $\max V
= \omega$. For each $v \in V$ such that $r(v)=1$ choose $\tilde{v}
\in v^+$ and define $f:\max V \rightarrow \omega$ as follows: if $v
\in w^+$ for $w \in V$, $r(w)=1$ then $f(v)=\tilde{w}$.

We will prove that the function $f$ fulfils the claim. Clearly, the
function $f$ is not constant on any $U \in u$. Consider $T= \{ v \in
V : r(v)=\omega \}$. It is sufficient to prove that
$r(f(v^\uparrow)) = \omega$ for each $v \in T$. Let $v^\uparrow =
(n)\looparrowleft(v_n)^\uparrow$ for $v \in T$. We have $r(\int
(v_n)^\uparrow) = r(f(\int (v_n)^\uparrow))+1$, so $\lim\limits_{n <
\omega}r(f(\int (v_n)^\uparrow)) =\omega = \lim\limits_{n <
\omega}r(\int (v_n)^\uparrow)$, and so $r(f(\int V)) = r(\int V) =
\alpha-1$.

Suppose that $f\mid_U$ is finite-to-one for some $U \in u$. Then
$\omega \setminus U \in \int V$, contradiction with $\int V \subset u$.

On the other hand by Theorem 2.3, $f(u) \in \cP_\gamma$ for a
certain $\gamma \leq \alpha$.  $\Box$

\begin{thm}\num
If $\alpha$ is a countably infinite successor ordinal and $u \in
\cP_\alpha$, then there exists a function $f:\omega \rightarrow
\omega$ such that:

1) $f^n\mid_U$ is not finite-to-one for any $n \in \omega$ and any
$U \in f^{n-1}(u)$ ($f^0(u)=u$),

2) the sequence $(f^n)_{n<\omega}$ is (pointwise) convergent;

3) $f^n(u) \in \cP_\alpha$ for each $n<\omega$, and
$\left(\lim_{n<\omega}f^n\right) (u)\in \cP_\alpha$.
\end{thm}

\Proof: Let ${V}$ be a monotone sequential cascade
 of rank $\alpha -1$ such that $\int {V} \subset u$.

Let $T= \{t \in  V : r(v)=\omega\}$. 
Without loss of generality we may assume that for each $v\in T$, for all $n<\omega$ each branch
of $v^\uparrow_n$ has length $r(v_n)$. For each $v\in V$ take a non decreasing sequence $(a^v_n)_{n<\omega}$ of natural numbers, such that $a^v_n \leq r(v_n)$, $\lim_{n\to\infty}a^v_n=\omega$, $\lim_{n\to\infty}(r(v_n)-a^v_n)=\omega$.

For each pair $(v,n)$ where $V\in T$, $n<\omega$ take a set $T_{v,n}=\{t\in v^\uparrow_n: r(t)=a^v_n\}$.
For each $t\in T_{v,n}$ take a function $f_{v,n}:\max t^\uparrow \to \max t^\uparrow$ defined like in the proof of case $4\Rightarrow 1$ in Theorem 3.1, and glue all this functions in a function $f:\max V \to \max V$ which satisfies a claim. $\Box$

\section{Relatively RK-$\alpha$-minimal points.}

Recall that a free ultrafilter $u \in \beta\omega $ is RK-minimal $^{1)}$ \footnotetext{$^{1)}$ also known as Ramsey ultrafilters or selective ultrafilters.} if for
each $f:\omega \rightarrow \omega$  either $f(u)$ is a principal
filter or $f(u)\approx u$. The existence of RK-minimal points is
%ZFC-independent
independent of ZFC
(see \cite{ComNeg}, \cite{vMill}). The following
theorem describes some properties of RK-minimal points.

\begin{thm} (see \cite{ComNeg}, \cite{vMill})
The following statements are equivalent for a free ultrafilter $u$ on $\omega$.

1) $u$  is RK-minimal;

2) For each function $f:\omega\rightarrow\omega$ there exists a set
$U \in u$ such that either $f\mid_U$ is constant or $f\mid_U$ is
one-to-one;

3) $u$ is a P-point and for each finite-to-one
function $f:\omega\rightarrow\omega$ there exists a set $U \in u$
such that $f\mid_U$ is one-to-one;

4) For each partition $d=\{d_n; n<\omega\}$ either there exists a
set $U \in u$ such that $\card(U \cap {d_n}) \leq 1$ for each $n <
\omega$, or there exists $n_0$ such that $d_{n_0} \in u$.
\end{thm}

An ultrafilter $u \in \cP_{\alpha}$ is referred to as {\it
relatively RK-$\alpha$-minimal} if for each $f:\omega \rightarrow
\omega$ there is either $u \approx f(u)$ or $f(u)\in \cP_\beta$ for
a certain $\beta < \alpha$; {\it relatively
$<_\infty$-$\alpha$-minimal} ultrafilters are those $u \in \cP_\alpha$
which each not finite-to-one (on each set $U\in u$) image is not in $\cP_\alpha$. 

The following two propositions, admitting straightforward proofs,
are useful in investigations of images of contours.

\begin{propo}\num
If $(p_n)$ is a sequence of filters,
%if
$p$ is a filter
%,
and
%if
$f:\omega \rightarrow \omega$
is a function, then $f(\int_pp_n) = \int_{F(p)}o_m$,
where $(o_n)$ is a sequence (possibly a finite sequence) of filters such that
$o_i \not= o_j$ for $i \not=j$ and $\{ o_j:j<\omega\} = \{ f(p_n) :
n<\omega \}$, $F(n)=i$ iff $f(p_n) = o_i$.
\end{propo}

Notice that $F$ depends on the order on the set $f(\{p_n : n<\omega
\})$, so in the remainder of this paper a function $F$ for $f$ is an
arbitrary (but fixed) function among such functions.

\begin{propo}\num
Let $(p_n)$, $p$, $(o_n)$ and $F$ be as in Proposition 4.2.
%If
Suppose that
there
exists a set $P \in F(p)$ such that the sequence $(o_i)_{i \in P}$
is discrete and
%if
there exists a set $H \in p$ such that $F\mid_H$ is one-to-one and
$p_n \approx o_{F(n)}$ for each $n \in H$. Then $\int_p p_n \approx
\int_{F(p)}o_i = f(\int_p p_n)$.
\end{propo}

\begin{thm}\num
Let $m < \omega$. If $(p_n)$ is a discrete sequence of relatively
RK-$m$-minimal free ultrafilters on $\omega$ and
%if
$p$ is an
RK-minimal free ultrafilter, then $\int_pp_n$ is relatively
RK-$m+1$-minimal.
\end{thm}

\Proof:  Let $p$ and $(p_n)$ be as in the assumptions.
Let $f$ be a
function $f:\omega \rightarrow \omega$. By Theorem 2.4 $\int_pp_n
\in \cP_{m+1}$. Take $(o_l)$ and $F$
%like
as
in Proposition 4.2. Thus,
$f(\int _pp_n) = \int_{F(p)} o_l$. Without loss of generality we
may assume that $\int_{F(p)} o_l \in \cP_{m+1}$.
We want to prove that $\int_{F(p)}o_i\not\in \cP_{m+1}$ or $\int_{F(p)}o_i
\approx \int_pp_n$. For this and consider
two
cases:
%possibilities $1$) $F(p)$ is a principal filter, or
%$2$) $F(p)$ is a free filter.

{\sl Case 1\/}.
%In case 1)
{\it $F(p)$ is a principal filter}.
In this case
there exists $i<\omega$ such that $\{i\}\in F(p)$ and
thus $o_i = \int_{F(p)}o_i$. Since $o_i = f(p_j)$ for some
$j<\omega$, and by Theorem 2.3 $o_i \in \cP_\beta$ for some $\beta
\leq m$,
%so
we have
$f(\int_pp_n) \not\in \cP_{m+1}$.

{\sl Case 2\/}.
{\it $F(p)$ is a free filter}.
%In case 2)
Then
$F(p)$ is a free ultrafilter, and thus, $F(p)\approx p$,
since $p$ is
RK-minimal.
%) $F(p)\approx p$.
Define sets $D_i = \{ n<\omega : o_n
\in \cP_i\}$, for $i\in \{1,...,m\}$. Since $o_n \in \bigcup_{j=1}^m
\cP_j$, there exists exactly one $i_0\in \{1,...,m\}$ such that
$D_{i_0} \in F(p)$.
%There are two subcases: either 2.1) $i_0 < m$ or
%2.2) $i_0=m$.

{\sl Subcase 2.1\/}.
$i_0 < m$.
%In case 2.1) let
Let us take a discrete sequence $(q_i)$ of ultrafilters
such that $q_i \approx o_i$, in this aim consider a partition
$(A_i)_{i<\omega}$ of $\omega$ into infinite sets, and a sequence $(f_n)$ of bijections
$f_i:\omega \to A_i$ and put $q_i=f_i(o_i)$. By Theorem 2.4 $ \int_{F(p)} q_i \in
\cP _{i_0+1}$ there is $ \int_{F(p)}o_i \preceq \int_{F(p)}q_n$
and so $\int_{F(p)}o_i \not\in
 \cP_{m+1}$.

%Thus we are left with case 2.2.
{\sl Subcase 2.2\/}.
$i_0=m$.
For each $i \in D_m$ and for each $n$ with $F(n)=i$ we have $p_n\approx f(p_n)=o_{F(n)}=o_i$
by RK-minimality of $p_n$.

Let $i_1=\min  D_m $. There exists a set $A_1$ such that $A_1
\in o_{i_1}$ and $A_1 \not\in \int_{F(p)} o_i$ (because we are not
in case 1). If numbers $i_r$ and sets $A_{r}$ for $r<t$ are
already defined, we define $i_t = \min \{i \in D_m:
(\bigcup_{r=1}^{t-1}A_{r})^c \in o_{i}\}$, and let $A_{t}$ be
such a set that $A_{t}\subset ( \bigcup_{r=1}^{t-1}A_{r} ) ^c$,
$A_{t} \in o_{i_t}$, $(A_{t})^c \in \int_{F(p)}o_i$ (such a set
exists because we are not in case 1, and $
(\bigcup_{r=1}^{t-1}A_{r})^c \in o_{i_t}$).
In this way we obtain a
sequence $(A_{r})_{r<\omega}$ of pairwise disjoint sets such that
$(A_{r})^c \in \int_{F(p)}o_i$ for each $r<\omega$, and for each
$i<\omega$ there exists a number $r<\omega$ such that $A_{r}
\in o_i$. Thus, the sequence $(A_{r})_{r<\omega}$ defines a
partition $s=(S_n)_{n<\omega}$ of $D_m$ into non-empty sets by
letting $i \in S_n$ if and only if $A_{n} \in 
o_i$. There is no $n$ such that
$S_n \in F(p)$ and $F(p)$ is RK-minimal, so by Theorem 4.1 there
exists a set $P\in F(p)$ with $P \subset D_m$ such that $\card(
P\cap S_n) \leq 1$ for each $n \in \omega$ (the sequence $(o_i)_{i
\in P}$ is discrete). The same Theorem 4.1 shows that there exists a
set $H \in p$ such that $F\mid_H$ is one-to-one. Without loss of generality
$F(H) \subset D_m$ and since $P_i \approx o_{F(i)}$ for all $i \in H$
we are in the assumption of Proposition 4.3 so we conclude
$\int_{F(p)}o_i \approx \int_pp_n$. $\Box$

By induction, by Theorem  4.4 one can easily prove the following
Corollary 4.5:

\begin{cor}\num
If there exist RK-minimal ultrafilters in $\beta\omega \setminus
\omega$, then for each $n<\omega$ there exist relatively RK-n-minimal
ultrafilters.
\end{cor}

In contrast to the above Corollary 4.5, for infinite $\alpha$'s
%, by
%Proposition 3.5
we obtain the following
from Theorem 3.6:
\begin{cor}
There are no $<_\infty$ (and so no RK)  relatively minimal
ultrafilters in classes of infinite successor index of the
P-hierarchy.
\end{cor}

\begin{prob}\num Do relatively
RK-$\alpha$-minimal ultrafilters exist for
limit ordinals $\alpha \leq \omega_1$? %ultrafilter in
\end{prob}

We may also consider RK-minimal elements in classes of ordinal
ultrafilters. An ultrafilter $u \in J^*_{\omega^{\alpha}}$ is
referred to as a {\it relatively ordinal RK-$\alpha$-minimal} if for
each $f:\omega \rightarrow \omega$ either $u \approx f(u)$
or $f(u)\in J^*_{\omega^{\beta}}$ for a certain $\beta < \alpha$.

One can get a very similar result to Theorem 4.4 for ordinal
ultrafilters:

\begin{thm}
Let $m < \omega$. If $(p_n)$ is a discrete sequence of relatively
ordinal RK-$m$-minimal free ultrafilters on $\omega$ and
%if
$p$ is a
RK-minimal free ultrafilter, then $\int_pp_n$ is relatively ordinal
RK-$m+1$-minimal.
\end{thm}

%A
The
proof
%, that
is very similar to the proof of Theorem 4.4
and
%,
uses the
following reformulation of a theorem of Baumgartner \cite[Theorem
4.2]{Baum}:

\begin{thm} Let $(\alpha_n)_{n<\omega}$ be a non-decreasing sequence of ordinals
less than $\omega_1$, let $\alpha = \lim\limits_{n < \omega}
(\alpha_n )$ and let $(X_n)$ be a partition of $\omega$. If $(p_n)$
is a sequence of ultrafilters such that  $X_n \in p_n \in
J^*_{\omega^{\alpha_n}}$ and $p$ is a P-point, then $\int_p p_n \in
J^*_{\omega^{\alpha+1}}$.
\end{thm}

Notice also that, by Proposition 3.2, each element of a class of
finite index of the P-hierarchy is relatively $<_\infty$-minimal.

In \cite[Theorem 3.3]{Laf} Laflamme  built (under MA for
$\sigma$-centered partial orderings) a special ultrafilter
%(call it
$u_0 \in J^*_{\omega^{\omega+1}}$
%),
the only RK-predecessor of which is a Ramsey ultrafilter. In
the proof of \cite[Theorem 3.13]{Star-P-hier} it is shown that $u_0
\in \cP_3$ and that $U_0$ is not in the form of a contour, note also that Laflamme's ultrafilter
$u_0$ is different then ultrafilter build in Theorem 4.4 which is a contour. Therefore,
we have:

\begin{thm}\num It is consistent with ZFC that there exists
a relatively-RK-3-minimal ultrafilter that is not in the form of a contour.
\end{thm}

\section{ Generic existence}

In this section, the P-hierarchy is understood as
$\bigcup_{1<\alpha<{\omega_1}} \cP_\alpha$.

Let $V$ be a cascade. Denote $r_\alpha(V) = \{v \in V : r(v) =
\alpha\}$. Let $h$ be a function with the domain $V$
%, and
 such that
$h(v) \in \omega$ for $v \in \max V$, and $h(v)$ is a filter on the
set $\omega$, otherwise. We define $\int^h V$ inductively as follows:
$\int^hv^\uparrow$ is a principal ultrafilter generated by $h(v)$, for $v \in \max V$.
If $\int ^hw^\uparrow$ is defined for all $v_n\in v^+$ then $\int ^hv^\uparrow = \int_{h(v)}\int v_n^\uparrow$.

Recall that the family $\mathcal{F}$ of functions $\omega
\rightarrow \omega$ is a {\it dominating family} if for each
function $g:\omega \rightarrow \omega$ there exists $f \in
\mathcal{F}$ such that $f(n)\geq g(n)$ for almost all $n<\omega$,
i.e. there is $n_0$ such that $f(n)\geq g(n)$ for all $n>n_0$.
The {\it dominating number} $\mathfrak{d}$ is the minimum of
cardinalities of dominating families, and $\mathfrak{c}$ is the
cardinality of the continuum.

We say that filters belonging to the family $\mathbb{F}$ {\it exist
generically} if each filterbase of size less than $\mathfrak{c}$ can
be extended to a filter belonging to $\mathbb{F}$.

In \cite{Brendle} Brendle showed that:

\begin{thm} \cite[part of Theorem E]{Brendle}
The following are equivalent:

(a) $\mathfrak{d}=\mathfrak{c}$;

(b) ordinal ultrafilters  exist generically.
\end{thm}

We obtain the same result for the P-hierarchy. For this we need to
prove the following theorem.

\begin{thm}
 For each ordinal $1<\alpha<\omega_1$
and for each monotone sequential contour $\cal V$ of rank $\alpha$, the
minimum of cardinalities of filterbases of $\cal V$ is $\mathfrak{d}$.
\end{thm}
\Proof: First, we will show that there exists a base of cardinality
$\mathfrak{d}$.  Let ${\cal V}= \int V$ for a monotone sequential cascade
$V$. Let $\mathbb{D}=\{ d_\beta: \beta < \mathfrak{d} \}$ be a
dominating family for functions $V \rightarrow \omega$ (there exists
such a family since $V$ is countable). Define family $\mathbb{F}=\{
f_{d,n}: d \in \mathbb{D}, n<\omega \}$ as follows: $f_{d,n}(v)=
d(v)$ for $v \not=\emptyset_V$ and $f_{d,n}(v)= n$ for $v
=\emptyset_V$. For each $v \in \max V$ let $(\emptyset_V=v^{0,v}
\sqsubseteq v_{1,v} \sqsubseteq \cdot \cdot \cdot \sqsubseteq
v_{n(v),v}=v)$ be a branch of $V$ with maximal node $v$.
For each function $f:V
\rightarrow \omega$ define sets $V(f)\subset\max V$ by the
condition: $v\in V(f)$ if and only if for each $i$
there exists $k$ such that $v^{i+1,v} = v^{i,v}_k)$ and $k
\geq f(v^{i,v})$. A typical base of the cascade $V$ is as follows:
$\{ V(g): g: V\rightarrow \omega \}$.
Take any $g:V\rightarrow \omega$. Since $ \mathbb{D}$ is a
dominating family, there exists $d_{\beta_0} \in \mathbb{D}$ such
that $g \leq^* d_{\beta_0}$. Thus, the set $A= \{v\in V: g(v) > d_{\beta_0}(v)\}$ is
finite, so we can define $n_0= \max \{n : A \cap
(\emptyset_V)_n)^\uparrow \not= \emptyset\} + 1$. Therefore,
$V(f_{d_{\beta_0},n_0}) \subset V(g)$, thus $\{V(f): F \in \mathbb{F}\}$is a base.

Now, let us assume that there exists a base $\mathbb{B}=\{B_\beta:
\beta < \gamma\}$ of $\cal V$ with $\gamma < \mathfrak{d}$. Since
$\{V(f): f: V \rightarrow \omega \}$ constitutes the base for
$\cal V$,
%then
for each $\beta < \gamma$ there exists $f_\beta$ such that
$f_\beta: V \rightarrow \omega$ and $W(f_\beta) \subset B_\beta$. Let
$\mathbb{G}=\{f_\beta: \beta < \gamma\}$. Since $\card(\mathbb{G}) <
\mathfrak{d}$,
% then
for each $v \in V$ such that $r(v)=2$ the family $\{f_\beta
\mid_{v^+}: \beta < \gamma \}$ is not a dominating family on the set
$v^+$. For each $v$ such that $r(v)=2$,
% let
%us
take a function $g_v:v^+ \rightarrow \omega$ such that $g_v
\not\leq^* f_\beta\mid_{v^+}$ for each $\beta < \gamma$. Now, let
$g: V \rightarrow \omega$ be a function that
$g(v)=g_{\tilde{v}}(v)$ if $r(v)=1$ and $v\in \tilde{v}^+$;
otherwise, $g(v) = 1$. We have $V(g)\in \int V = {\cal V}$ and $V(g)
\not\supset V(f_\beta)$ for each $\beta<\gamma$. \Box

The {\it supercontour} is a filter of type
$\bigcup_{\alpha<\omega_1} {\cal V}_\alpha$, where
$({\cal V}_\alpha)_{\alpha<\omega_1}$ is an increasing sequence of monotone
sequential contours such that $r({\cal V}_\alpha)=\alpha$.

\begin{cor} Generic existence of the P-hierarchy is equivalent to
$\mathfrak{d}=\mathfrak{c}$.
\end{cor}
\Proof: By the Ketonen's theorem \cite[Theorem 1.1]{Ket} generic existence of
P-points is equivalent to $\mathfrak{d}=\mathfrak{c}$. By
Proposition 2.2 P-points belong to $\cP_2$ class; thus, if
$\mathfrak{d}=\mathfrak{c}$, then the P-hierarchy generically exists.

Let $\mathfrak{d} < \mathfrak{c}$, and
%let us
take an increasing
$\omega_1$-sequence $({\cal V}_\alpha)$ of monotone sequential contours such
that $r({\cal V}_\alpha)=\alpha$. Let $\mathbb{B}_\alpha$ be a base of
cardinality $\mathfrak{d}$ of ${\cal V}_\alpha$ (there exist by
Theorem 5.2). Let
$\mathbb{B}=\bigcup_{\alpha<\omega_1}\mathbb{B}_\alpha$.
Obviously,
$\card(\mathbb{B})=\mathfrak{d}$ and $\mathbb{B}$ is the base for a
supercontour, so it cannot be extended to any element of the
P-hierarchy. \Box

By Ketonen's Theorem \cite{Ket} each ultrafilterbase of cardinality
less than $\mathfrak{d}$ is the base of a P-point. In order to
obtain a similar
result
for other classes, we need the extra
assumptions that such a base can be extended to infinitely many
ultrafilters. To prove
this
%the result
 we need to quote the following
two results:

\begin{thm}\cite[Theorem 4.1]{Baum}
The $J^*_{\omega^2}$-ultrafilters are the P-point ultrafilters.
\end{thm}

We say that families $u$ and $o$ {\it mesh} (and we write $u \# o$)
whenever $U\cap O \not= \emptyset$ for every $U \in u$ and $O \in
o$.

\begin{propo} The following statements are equivalent:

a) For each successor ordinal $1<\alpha<\omega_1$ each filterbase of
cardinality less than $\mathfrak{c}$ which can be extended to
infinitely many ultrafilters, can also be extended to some elements of
$\cP_\alpha$;

b)  For each successor ordinal $1<\alpha<\omega_1$ each filterbase
of cardinality less than $\mathfrak{c}$ which can be extended to
infinitely many ultrafilters, can also be extended to some elements of
$J^*_{\omega^\alpha}$;

c) $\mathfrak{d}=\mathfrak{c}$.
\end{propo}

\Proof: For $\mathfrak{d} < \mathfrak{c}$ a proof is analogical to the second part of the proof of Theorem 5.3,
with an additional use of Proposition 2.8 for the case of ordinal ultrafilters.

Now let $\mathfrak{d}=\mathfrak{c}$, and let $\mathbb{B}$ be a
proper filterbase of cardinality $< \mathfrak{c}$. By the Ketonen's Theorem
\cite{Ket} $\mathbb{B}$ can be extended to a P-point, and so we can
assume that $\alpha > 2$. 
by the assumption, there exists a family $\{F_n\}_{n<\omega}$
of pairwise disjoint sets such that $F_n \# \mathbb{B}$ for each
$n<\omega$. Let $(p_n)$ be a sequence of P-points such that $\mathbb{B}\cup \{F_n\} \subset p_n$. Take a monotone sequential cascade $V$ of rank $\alpha$. Put $R=\{v\in V: r(v)=1\}$ and without loss of generality assume that for each $v \not\in R$ a cascade $v^\uparrow$ has no branches of length $1$. 
Let $g$ be an arbitrary bijections $g: R \to \omega$ and let $f_v: v^+\to F_n$ be an arbitrary bijection for each $v \in R$.
Let $h$ be a function which domain is $V$, defined as follows:

$h(v')=f_v(v')$ for $v'\in v^+$, $v \in R$

$h(v)=p_{g(v)}$ for $v\in R$

$h(v)=p_1$ for other $v\in V$.

Consider $\int ^hV$ and note following facts:

1) $\int ^hV \# \mathbb{B}$, since $p_n\# \mathbb{B}$ for all $n$;

2) $\int ^hV \in \cP_\alpha$, inductively by Theorem 2,5;

3) $\int ^hV \in J^*_{\omega^{\alpha+1}}$ by Corollary 5.4 and inductively by Theorem 4.8. $\Box$

\section{Existence}

We say that a cascade $V$ is {\it built by destruction of nodes of rank 1} in a cascade $W$ of rank $r(W)\geq2$
 iff for a set $R=\{w \in W: w(w)=1, r(w^-)=2\}$ there is: $V = W \setminus R$ and if $v\in R^{-W}$ then
$v^{+V} = (v^{+W}\setminus R) \cup (v^{+W}\setminus R^+)$, i.e. order on the cascade is unchanged.

Observe that if $W$ is a monotone sequential cascade then $V$ is also a monotone sequential cascade and if $r(W)$ is finite then $r(V)=r(W)-1$, if $r(W)$ is infinite, then $r(V)=r(W)$.

Assume that we are given a cascade of rank $\alpha$ and an ordinal $1<\beta \leq\alpha$. We shall describe an operation of {\it decreasing the rank} of a cascade $W$. 
The construction is inductive: 

For finite $\alpha$, we can decrease rank of $W$ from $\alpha$ to $\beta$ by applying $\alpha - \beta$ times an operation of destroying  nods of rank 1 (u.e. if $\alpha=\beta$ then the cascade is unchanged).

For infinite $\alpha$. Suppose that for each pair $(\delta, \gamma)$ where $1<\delta\leq\gamma<\alpha$, and for each cascade $W$ of rank $\gamma$ the operation of decreasing of the rank of $W$ from $\gamma$ to $\delta$ is defined. Let $W$ be a monotone sequential cascade of rank $\alpha$, let $(\beta_n)$ be a nondecreasing sequence of ordinals such that: 
$\beta_n=0$ if and only if $r(W_n)=0$, $\beta_n \leq r(W_n)$ and $lim_{n\to \infty}(\beta_n+1)=\delta$. Let, for each $n<\omega$, $V_n$ be the cascade obtained by decreasing of rank of $W_n$ to $\beta_n$. Finally let $V=(n) \looparrowleft V_n$.

Clearly for infinite $\alpha$ the operation of decreasing of rank is not defined uniquely. 
Observe also that the above described decreasing of rank of a cascade $W$ does not change $\max W$.
If a cascade $V$ is obtained from $W$ by decreasing of rank, then we write $V \triangleleft W$. 
Trivially $V \triangleleft W$ and inductively  $\int V \subset \int W$.
\vspace{0.5cm}

\begin{thm} \cite{Dol-multi} 
If $({\cal V}_n)_{n<\omega}$ is a sequence of monotone sequential contours of rank less than $\alpha$ and 
$\bigcup_{n<\omega} {\cal V}_n$ has the finite intersection property, then there is no monotone sequential contour ${\cal W}$ of rank $\alpha +1$  such that ${\cal W} \subset \langle \bigcup_{n<\omega} {\cal V}_n \rangle $.
\end{thm}

Before we prove the main Lemma we shall prove a technical claim;

\begin{lemm}
Let $V$ be the cascade of rank $\alpha$, $W$ be cascade obtained from $V$ by decreasing the rank of $V$ to $\beta <\alpha$ and let $\beta< \gamma < \alpha$. Then there is a cascade $T$ of rank $\gamma$ such that $W\triangleleft  T \triangleleft  V$.  
\end{lemm}

\Proof: If $\beta=1$ then it suffice to take any monotone sequential cascade $T$ obtained by decreasing of the rank of $W$ to $\gamma$.

If $\beta>1$ then take $(\beta_n)$  - a nondecreasing sequence of ordinals such that: 
$\beta_n=0$ if and only if $r(W_n)=0$, $r(V_n) \leq \beta_n \leq r(W_n)$ and $lim_{n\to \infty}(\beta_n+1)=\gamma$.
By inductive assumption one can find $(T_n)$ a sequence of monotone sequential contours such that $V_n \triangleleft T_n \triangleleft W_n$. Put $T=(n) \looparrowleft T_n$. $\Box$

We write $V \blacktriangleleft_1  W$ if $\max W \in \int V$ and $V^{\downarrow \max W} \triangleleft W$.
We write $V \blacktriangleleft_2  W$ if $\max V \in \int W$ and $V \triangleleft W^{\downarrow \max V}$.  Trivially Lemma 6.2 is true also for $\blacktriangleleft_1$, $\blacktriangleleft_2$ instead of $\triangleleft$. 

 \begin{lemm}\num
Let $\alpha < \omega_1$ be a limit ordinal and let $({\cal V}^n : n<\omega)$ be
a sequence of monotone
sequential contours such that $r({\cal V}^n)<r({\cal
V}^{n+1})<\alpha$ for every $n$ and such that $\bigcup_{n < \omega}{\cal V}^n$ has the finite intersection property. 
Then there is no monotone sequential contour
$\cal{W}$ of rank $\alpha$ such that ${\cal{W}} \subset
\langle \bigcup_{n<\omega} {\cal V}^n \rangle$.
\end{lemm}

\Proof: Put $\alpha_n=r({\cal V}^n)$, without loss of generality we may assume that $\alpha_1\geq 3$. Assume that there exists a monotone sequential contour ${\cal W}$ of rank $\alpha$ such that 
${\cal W} \subset \langle \bigcup_{n<\omega} {\cal V}^n \rangle$. 
We build a cascade $W$ and a sequence of cascades $(W^n)_{n<\omega}$ such that:
\begin{itemize}
\item $\int W ={\cal W}$;
\item $W^n \blacktriangleleft_1 W^{n+1}$ for all $n$;
\item $W^n \blacktriangleleft_2 W$ for all $n$;
\item $r(W^n) =  \alpha_n +3$ for all $n$;
\item $r(W^n_i) = \alpha_n+2$ for all $n$ and all $i$;
\item $r(W^n_{i,j}) = \alpha_n+1$ for all $n$, $i$ and $j$..
\end{itemize}
\vspace{0.5cm}

Fix any monotone sequential cascade $\bar{W}$ such that $\int \bar{W} = {\cal W}$.
Let $\bar{W}^m$ be the cascade obtained from $\bar{W}$ by cutting every subcascade $\bar{W}_i$ of rank smaller than $\alpha_m +2$ and every subcascade 
$\bar{W}_{i,j}$ of rank smaller than $\alpha_m +1$. Observe that we cut only finitely many subcascades $\bar{W}_i$ and for the other $\bar{W}_i$ only finitely many subcascades 
$\bar{W}_{i,j}$. Thus $\int \bar{W}^m = \int \bar{W} = {\cal W}$ for every $m$.   
\vspace{0.5cm}

Let $W=\bar{W}^1$ and $W_1$ be a cascade obtained from $\bar{W}_1$ by decreasing ranks of $W^1_{i,j} $  to $\alpha_1 +1$. Thus 
$W^1 \blacktriangleleft_2 W $. 

Since cascades $\bar{W^n}$ and $W^n$ are subcascades of $W$ thus for nodes (and so subcascades) of $\bar{W^n}$ and $W^n$ we may keep the indexation from $W$, to avoid the collision of notation we put those indexes in parenthesis.

Assume that $W^1 \blacktriangleleft W^2 \blacktriangleleft \ldots \blacktriangleleft  W^m $ have been defined. % such that $W_l \blacktriangleleft W$ and $r(W_l^{(i)(j)})= \alpha_l +1$ (thus $r(W_l)= \alpha_l +3$)   for every $l\leq m$.
We apply Lemma 6.2 to cascades    $W^m_{(i,j)}$ and $\bar{W}^{m+1}_{(i,j)}$ to define $W^{m+1}_{(i,j)}$ of rank $\alpha_{m+1}+1$ for those $(i,j)$ that $w_{i,j}\in W^m \cap \bar{W}^{m+1}$.

Let $K^{m+1}$ be a subcascade of $W$ with $K^{m+1}=\{\emptyset_W\}\cup  (\emptyset_W)^{+\bar{W}^{m+1}} \cup ((\emptyset_W)^{+\bar{W}^{m+1}})^{+\bar{W}^{m+1}}$.
Put $W^{m+1}=K^{m+1} \looparrowleft W^{m+1}_{(i,j)}$.

\vspace{0.5cm}

Next we build a decreasing sequence $(U_n)_{n<\omega}$ satisfying conditions $U_A$-$U_D$:
\begin{enumerate}
\item $U_A(n)$: $U_n \in \int W^n$; 
\item $U_B(n)$: $U_n \notin \langle \bigcup_{i\leq n} {\cal V}^i  \rangle $; 
\item $U_C(n)$: $U_n \cap (\omega \setminus \max {W}^{n+1}) = U_{n+1} \cap (\omega \setminus \max \bar{W}^{n+1})$;
\item $U_D(n)$: $U_n  \cap \max W_i \in \int W_i$ for all $n$ and all $i$. 
\end{enumerate}

In this aim first we built an additional sequence $(\widetilde{W}^n)$ of cascades by $\widetilde{W}^n = W^n \setminus \emptyset_{W^n}^+$ such that $ \emptyset_{\widetilde{W}^n}^{+} = \bigcup \left\{ w^{+} : w \in  \emptyset_{W^n}^{+}    \right\}.$ and that the rest of cascades we leave unchanged (we may say that $\widetilde{W}^n$ is obtained from $W^n$ bay destroying all nodes of rank $\alpha_m+2$). Notice that $\widetilde{W}^n$ is a monotone sequential cascade of rank $\alpha_m+2$, and that if a set $U_n$ fulfills conditions $U_B(n)$, $U_C(n)$ and belongs to $\widetilde{W}^n$ then the same set $U_n$ fulfills all conditions $U_A(n)-U_D(n)$.

\vspace{0.5cm}

Put $U_0 = \omega$. Assume that $U_0, U_1, \ldots , U_{n-1}$ was defined, but it is impossible to define $U_n$. This means that every set $U\in \int \widetilde{W}_n $ is contained in  $\langle \bigcup_{i< n } {\cal V}_i \rangle$. 
On the other side $\max \widetilde{W}_n \in {\cal W}$ and so the family
$\{ U\cap  \max \widetilde{W}_n : U \in \bigcup_{i\leq n} {\cal V}_i \} $ has the finite intersection property. By the theorem of Dolecki $\langle \{ U\cap  \max \widetilde{W}_n : U \in \bigcup_{i\leq n} {\cal V}_i \}\rangle $ do not contain any monotone sequential contour of rank $\alpha_n+2$ and so do not contain $\int \widetilde{W}_n$. A contradiction.
On each step of induction we can put $\bigcap_{i\leq n} U_i$ instead of $U_n$ and assume that the sequence $(U_n)_{n<\omega}$ is decreasing. %[NIE TAK SZYBKO! NIE KUMAM!] 
\vspace{0.5cm}

Notice that $\bigcup_{n<\omega} (\max {W}_{n+1})^c = \max W$, let $U= \bigcap_{n<\omega} U_n$. Conditions (1)-(4) guarantee that 

1) $U\in \int W$ and

2) $U \notin \langle \bigcup_{n<\omega} {\cal V}_n \rangle$.

To see 1) fix any $t<\omega$, note that $\max W^m \in \int W_t$ only for finite number of $m$. So the sequence $(U_n\cup R)_{n<\omega}$ is (decreasing and) almost constant on some $R\in \int W_t$. Therefore $\bigcap_{n<\omega}U_n \cap R$ is indeed a finite intersection of $R$ and $U_n$ all of which by condition $U_D$ belongs to $\int W_t$. So $\bigcap_{n<\omega}U_n \in \int W_t$ for all $t$, and so $U \in \int W$.

To see 2), assume that $U \in \langle \bigcup_{n<\omega} {\cal V}_n \rangle$, then there is a finite $M<\omega$ such that 
$U \in \langle \bigcup_{n<M} {\cal V}_n \rangle$. But $U_M  \notin \langle \bigcup_{n\leq M} {\cal V}_n \rangle$ and $U \subset U_M$. Thus $U \notin \langle \bigcup_{n<\omega} {\cal V}_n \rangle$. A contradiction. $\Box$

\begin{propo}\num\cite[part of corollary 2.6]{Star-P-hier} (ZFC)
Classes $\cP_1$ and $\cP_{\omega_1}$ are nonempty.
\end{propo}

\begin{thm}\num (CH)
Each class of the P-hierarchy is nonempty.
\end{thm}

\Proof: For successor $\alpha$'s and for 1 for $\omega_1$ we deal in Proposition 6.4. 
By well known result of W. Rudin CH implies existing of P-points so for successor $\alpha$ we are done
by Theorem 6.5.  Let $\alpha < \omega_1$ be limit ordinal. Let $({\cal
V}_n)$ be an increasing sequence of monotone sequential contours
such that $r({\cal V}_n)$ is an increasing sequence with
$\lim_{n<\omega} r({\cal V}_n) =  \alpha$. By CH we can order all
$\alpha$-partitions in an $\omega_1$ sequence $(P_\beta)$.

We will build a sequence $(Q_\beta)_{\beta<\omega_1}$ of subsets
of $\omega$ such that $Q_\beta$ is residual for the partition
$P_\beta$ and a family $\{Q_\beta: \beta < \omega_1\} \cup \bigcup_{n<\omega} {\cal V}_n$ has the finite intersection property.
Since $\bigcup_{n<\omega} {\cal V}_n$ is a filter and, by the Lemma 6.3 above, does not
contain any monotone sequential contour of rank $\alpha$,
thus there
exists a set $Q_1$ residual for the partition $P_1$ such that the
family $\{Q_1\} \cup \bigcup_{n<\omega} {\cal V}_n$ has the finite
intersection property. Suppose now that the sequence
$(Q_\beta)_{\beta<\gamma}$ is already built. If $\gamma<\omega$ then
consider the sequence $( {\cal V}_n \mid_{ \bigcap_{\beta<\gamma}
Q_\beta})_{n<\omega}$, this is an increasing sequence of monotone
sequential contours with $r({\cal V}_n) = r({\cal V}_n \mid_{\bigcap_{\beta<\gamma} Q_\beta})$
 thus by the Lemma 6.3 there exist
a set $Q_\gamma$ residual for the partition $P_\gamma$ and such that a
family $\{Q_\gamma\} \cup \bigcup_{{n<\omega} }  ({\cal V}_n \mid_{
\bigcap_{\beta<\gamma} Q_\beta})$ has the finite
intersection property and thus also a family $\{Q_\beta :
\beta\leq \gamma\} \cup \bigcup_{n<\omega} {\cal V}_n$ has the
finite intersection property. If $\gamma \geq \omega$ then we
enumerate the sequence $(Q_\beta)_{\beta<\gamma}$ by natural numbers
and obtain the sequence $(Q^{\gamma,n})_{n<\omega}$. Consider the
sequence $({\cal V}_n\mid_{ \bigcap_{m\leq n}
Q^{\gamma,m}})_{n<\omega}$, this is an increasing sequence of monotone
sequential contours with $r({\cal V}_n) = r({\cal V}_n\mid_{
\bigcap_{m\leq n} Q^{\gamma,m}})$. Thus by the Lemma 6.3 there exist a
set $Q_\gamma$ residual for the partition $P_\gamma$ and such that a
family $\{Q_\gamma\} \cup     \bigcup_{n<\omega} ({\cal V}_n\mid_{ \bigcap_{m\leq
n} Q^{\gamma,m}})$ has the finite intersection property
and thus also a family $\{Q_\beta: \beta\leq \gamma\} \cup \bigcup_{n<\omega}{\cal V}_n$ has the finite intersection property.
Thus a sequence $(Q_\beta)_{\beta<\omega_1}$ with described
properties exists.

 Now it is sufficient to take any
ultrafilter $u$ that contains $\{Q_\beta: {\beta<\omega_1}\} \cup
\bigcup_{n<\omega}{\cal V}_n$. Since $ \bigcup_{n<\omega}{\cal V}_n  \subset u$ then $u$
contains a monotone sequential
contour of each rank less then $\alpha$. Since $u$ contains
$\{Q_\beta: {\beta<\omega_1}\}$  thus $u$ contains residual set
for each $\alpha$-partition, and thus $u$ do not contain any
monotone sequential contour of rank $\alpha$. $\Box$

Notice that it was also shown

\begin{thm}\cite[reformulation of Theorem
3.12]{Star-P-hier}

MA$_{\sigma-{\rm center.}}$ implies $\cP_{\alpha+\omega}\not=
\emptyset$.

\end{thm}

It is worth to compare the above results with \cite[Theorem 4.2]{Baum},
where Baumgartner proved that
if P-points exist then for each successor $\alpha<\omega_1$
the class of $J^*_{\omega^\alpha}$ ultrafilters is nonempty, and with our theorem
from  \cite{Star-OCU} where we proved (in ZFC) that a class of $J^*_{\omega^\omega}$ ultrafilters is empty.

 \hspace{-2mm}
\gd

\medskip

{\small\sc \noindent {Wydzia\l}  {Matematyki Stosowanej},
Politechnika \'{S}l\c{a}ska, Gliwice, Poland

E-mail:  andrzej.starosolski@polsl.pl}

\end{document}